%% file: substitution-search-paper.tex
\documentclass[a4paper,12pt%
]{amsart}
\pdfoutput=1
\usepackage{amsmath, amssymb, amsthm, amscd}
\usepackage{graphicx}
\usepackage{enumerate}
\usepackage{caption}
\usepackage{subcaption}
\usepackage[all,cmtip]{xy}

\usepackage{algorithm}
\usepackage{algpseudocode}

\usepackage{etex} 

\usepackage[english]{babel}
\usepackage[latin1]{inputenc}
\usepackage[all,cmtip]{xy}
\usepackage{times}
\usepackage{cancel}

\theoremstyle{remark}
\newtheorem*{remark*}{Remark}

\numberwithin{equation}{section}
\numberwithin{theorem}{section}

\usepackage{macros}

\title[Substitution search]{A computer search for planar substitution tilings with $n$-fold rotational symmetry}

\author{Franz G\"{a}hler}
\address{Bielefeld University}
\thanks{The first and third authors were partly supported by the German Research Council (DFG), CRC 701.  The third author was also partly supported by the Fields Institute during a research visit.}

\author{Eugene E. Kwan}
\address{Harvard University}

\author{Gregory R. Maloney}
\address{Newcastle University}

\begin{document}

\date{\today}

\begin{abstract}
We describe a computer algorithm that searches for substitution rules on a set of triangles, the angles of which are all integer multiples of $\pi /n$.  
We find new substitution rules admitting $7$-fold rotational symmetry at many different inflation factors.  
\end{abstract}

\maketitle

\input{intro}
\input{lengths-areas}
\input{algorithm}
\input{results}

\bibliographystyle{abbrv}
\bibliography{bib-cohomology}
\end{document}

%% file: intro.tex
\section{Introduction}\label{SEC:intro}

\subsection{Motivation}\label{SUBSEC:motivation}

In \cite{DN:n-fold} a method is described for constructing substitutions on the set of all triangles with angles that are integer multiples of $\pi/n$, subject to an appropriate normalization.  
This method gives rise to an infinite family of planar substitution tilings, which furnish examples that are of interest for their own inherent structure, for their roles as models of aperiodic solids called quasicrystals, and for aesthetic reasons.  
But the tilings constructed via this method do not appear to exhibit global $n$-fold rotational symmetry (except in the special case $n=9$), although the underlying tiling space is $n$-fold symmetric.  

Still there is another substitution described in \cite[Figure 12]{DN:n-fold} that is defined on a proper subset of the triangles with angles that are integer multiples of $\pi/7$.  
This substitution does not arise from the general construction; indeed, the method of its discovery is not explained, yet it appears, upon casual inspection, to give rise to a tiling with global $7$-fold rotational symmetry.  
This turns out upon closer inspection to be false, as certain isosceles triangles appear in reflected positions, breaking the symmetry.  
The authors of \cite{DN:n-fold} also observe that this last substitution is special in that it admits a local matching rule (see \cite{GS:matching-rules} and \cite{TE}) whereas, in all of the cases that they checked, the substitutions arising from their general method do not.  

The goal of this work is to search for other substitutions that are similar to this extra substitution in that they are defined on a proper subset of the triangles with angles that are integer multiples of $\pi/n$.  
In particular, the intention is that, by selecting a minimal subset of these triangles, the substitutions found will produce at least one tiling possessing global $n$-fold rotational symmetry.  
Such substitutions will necessarily not arise from the general construction in \cite{DN:n-fold}, because the general construction uses all triangles with angles that are integer multiples of $\pi/n$, at least for $n$ not divisible by $3$.  

It is of particular interest to find multiple different substitution rules on the same set of prototiles that are compatible with one another, meaning that they can be combined to produce edge-to-edge tilings.  
There has been much recent work on different tilings spaces that arise from the combination of two or more substitution rules on the same prototile set.  
(Some terms mentioned here, such as ``prototile,'' will be defined formally in Section \ref{SUBSEC:definitions}.)  
These fall into two classes: the multi-substitution tilings (see \cite{F:s-adic}, \cite{D:lin-recurrent}, \cite{PFS:fusion}, \cite{GM:multi-one-d}, \cite{G:fivefold}, and \cite{PV:frequency}), that are obtained by choosing a substitution for each hierarchical level and applying it to all tiles at that level; and the random substitution tilings (see \cite{GL:random}, \cite{N:random-one}, \cite{N:random-two}, and \cite{BM:noble-means}), that are obtained by making separate choices of substitution for each tile at each hierarchical level.  
While it is easy to find examples of such families of substitutions in one dimension, in two dimensions it is harder.  
Most known examples are either constant length substitutions \cite{F:multi-dim-constant-length} or lack the edge-to-edge property \cite{GL:random}.  

The edge-to-edge property is obviously desirable from the point of view of modeling quasicrystals.  
But constant length substitutions have integer inflation factors, and so exhibit behavior markedly different from that of substitutions with non-integer inflation factors (see \cite{K:constant-length} and \cite{D:constant-length}).  
Therefore there is a need for a collection of examples in two dimensions that are edge-to-edge and have non-integer inflation factors.  
This project addresses that need.

\subsection{Background}\label{SUBSEC:background}

Other projects have been undertaken with a similar goal in mind.  
In \cite{F:thesis}, a family of substitution rules is introduced, one for each $n>7$, that generalises the extra substitution in \cite[Figure 12]{DN:n-fold}.  
This involves amalgamating adjacent triangles into quadrilaterals and pentagons to bypass a negative area obstruction.  
In \cite{H:rhomb}, a family of substitutions on rhombic tiles is introduced, generalising a rule of Goodman-Strauss to orders of symmetry greater than $7$.  
In the notation of Section \ref{SUBSEC:lengths}, the inflation factor in \cite{F:thesis} is $1+a_2$ and the inflation factor in \cite{H:rhomb} is $2+a_2$, where $a_2=2\cos(\pi/n)$.  
The former is a unit in its ring of integers, but not a Pisot-Vijayaraghavan (PV) number, except in a few cases.  
The latter is neither a PV number nor a unit in its ring of integers.  
Representative pictures of both families of rules can be found at \cite{TE} under the names ``cyclotomic trapezoids'' and ``Harriss's $9$-fold rhomb'' respectively.  

Both of these works succesfully adapt a substitution rule that was originally defined for $n=7$ to arbitrary $n$, therefore producing an infinite family of substitutions.  
In the process, both of them introduce extra prototiles, so that they no longer work with a minimal set.  
Also, the tilings that result from the substitution rules they describe do not exhibit local $n$-fold symmetry in the $k$-fold substituted image of any prototile, and hence the associated tiling spaces do not contain any tiling with global $n$-fold symmetry.  

The approach in this work is exploratory rather than constructive.  
In \cite{F:thesis} and \cite{H:rhomb} the approach was to generalise to higher orders of symmetry $n$ a substitution rule that has already been found at a low order of symmetry, thereby producing an infinite family of substitution rules, parametrized by $n$.  
Here we make no attempt to build on rules that already exist, but instead search exhaustively for all rules that fit a certain description.  
In this way we find many rules that do not appear in \cite{F:thesis} or \cite{H:rhomb} because they do not fall into those families.  
The drawback is that we can only ever hope to find finitely many substitution rules by this method---although, as we shall see, this finite list is quite long, and can easily be made longer.  

Since we impose no preconditions on the inflation factors that we search, accordingly we must search a large parameter space.  
To do this it is necessary to use a computer.  

\subsection{Definitions}\label{SUBSEC:definitions}

A \defemph{tile} is a subset of $\R^d$ homeomorphic to the closed unit disk.  
A \defemph{patch} is a collection of tiles, any two of which intersect only in their boundaries.  
The \defemph{support} of a patch is the union of the tiles that it contains.  
A \defemph{tiling} is a patch, the support of which is all of $\R^d$.  

Let us restrict our attention to the case $d = 2$, and let us consider only polygonal (in fact, triangular) tiles having some finite set $\PT$ of representatives up to isometry.  
The elements of $\PT$ are called \defemph{prototiles}, and we denote by $\PatchSet(\PT)$ the set of all patches consisting of tiles that are congruent to these prototiles.  
Let us also suppose that we have a \defemph{substitution}, that is, a map $\sub : \PT \to \PatchSet(\PT)$ for which there is an \defemph{inflation factor} $\infl > 1$ such that, for any $\prototile\in\PT$, the support of $\sub (\prototile)$ is $\infl \prototile$.  
Letting $\PT = \{ \prototile_1,\ldots , \prototile_k\}$, we define the \defemph{substitution matrix} of $\sub$ to be the $k \times k$ integer matrix, the entry of which at position $(i,j)$ is the number of tiles isometric to $\prototile_i$ in $\sub(\prototile_j)$.  

\begin{figure}
\includegraphics[width=\textwidth]{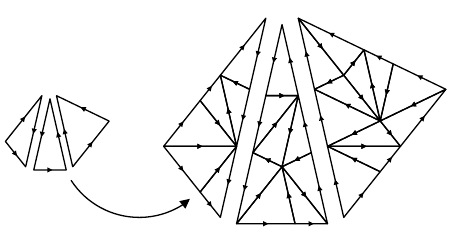}
\caption{A substitution on three prototiles, with edge orientations.\label{FIG:sample-substitution}}
\end{figure}

Figure \ref{FIG:sample-substitution} depicts a substitution on three prototiles, with the arrows on the edges of the tiles describing \defemph{edge orientations}, which are defined in Section \ref{SUBSEC:objective} below.  
This substitution has the same prototile shapes and substitution matrix---that is, 
\[
\left[ \begin{array}{ccc} 3&3&5\\ 1&4&3\\ 2&1&3 \end{array}\right]
\]
---as the substitution in \cite[Figure 12]{DN:n-fold}, although it is indeed a different substitution.  
In particular, the prototiles, although congruent to those of \cite[Figure 12]{DN:n-fold}, have different edge orientations, and so should be seen as different prototiles.  

\begin{figure}
\includegraphics[width=\textwidth]{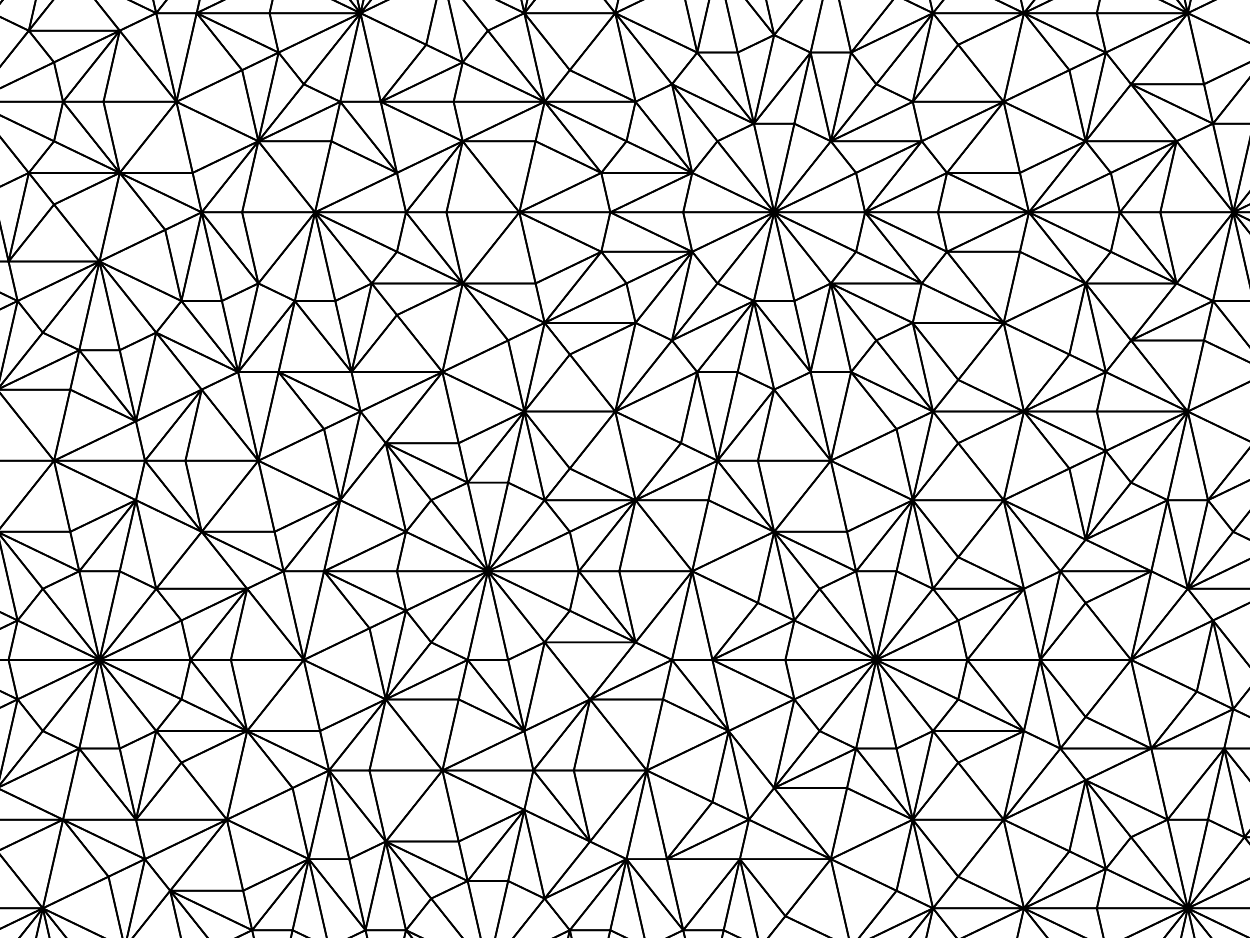}
\caption{A substitution tiling constructed from the substitution in Figure \ref{FIG:sample-substitution}.\label{FIG:sample-tiling}}
\end{figure}

We can extend the definition of $\sub$ to all tiles isometric to the elements of $\PT$ in the following way.  
If $\prototile$ is a prototile, $\trans$ is an orthogonal transformation of $\R^2$, and $\vectorv\in\R^2$ is a translation vector, then $\sub (\trans(\prototile)+\vectorv) = \{ \trans(\tile)+\infl \vectorv \ | \ \tile \in \sub(\prototile)\}$.  
Once we have done this, we can extend the definition of $\sub$ further to include all patches in $\PatchSet(\PT)$ by declaring $\sub (\Patch) = \{ \sub(\tile) \ | \ \tile \in \Patch\}$.  
The chief motivation for this is that, under certain easily-satisfied conditions, we can produce a tiling by starting with some prototile $\prototile$ and repeatedly applying $\sub$ \cite{GS:book}.  
Let us refer to such tilings as \defemph{substitution tilings}.  
Figure \ref{FIG:sample-tiling} depicts a substitution tiling constructed in this way from the substitution in Figure \ref{FIG:sample-substitution}.  
It clearly possesses local $7$-fold rotational symmetry.   

\subsection{The Objective}\label{SUBSEC:objective}

Our purpose here is to find substitutions that obey a certain rule and use certain prototile sets.  
The rule is that the resulting substitution tilings must be \defemph{edge-to-edge} \cite{GS:book}, which means that, if two tiles in a tiling intersect, their intersection is a face of both tiles.  
For polygonal tiles in two dimensions, this amounts to the property that, if the vertex of one tile touches another tile, then it is also a vertex of the other tile.  

The edge-to-edge condition for substitution tilings arising from $\sub$ is equivalent to the condition that the patches $\sub^m(\prototile)$ be edge-to-edge for all $m\in\N$ and all prototiles $\prototile$.  
Verifying this condition seems at first to require checking infinitely many patches, but in fact it is sufficient to check for all prototiles $\prototile$ that the patches $\sub(\prototile)$ are edge-to-edge and satisfy one additional condition, which is described below.  

Given a tile $\tile$ with an edge $\edge$ having end points $\vectorv_1$ and $\vectorv_2$, a \defemph{edge orientation} is a map $\eo_\tile$ that assigns to $\edge$ one of its vertices $\vectorv_i$.  
We can represent $\eo_\tile$ graphically by drawing an arrow on $\edge$ originating at $\eo_\tile(\edge)$ and terminating at the other vertex, as has been done in Figure \ref{FIG:sample-substitution}.  
Then we require that all prototiles $\prototile$ and tiles $\tile$ have edge orientations on all of their edges satisfying the following conditions.  

\begin{enumerate}
\item  (Isometry equivariance)  If $\tile = \trans(\prototile) + \vectorv$, then, for any edge $\edge$ of $\prototile$, the corresponding edge in $\tile$ has the same edge orientation; that is, $\eo_\tile(\trans(\edge)+\vectorv) = \trans(\eo_\prototile(\edge))+\vectorv$.  
\item  (Matching)  If an edge $\edge$ lies in two tiles $\tile_1$ and $\tile_2$, then it receives the same edge orientation from both of them; that is, $\eo_{\tile_1}(\edge) = \eo_{\tile_2}(\edge)$.  
\item  (Preservation under $\sub$)  If two prototiles $\prototile$ and $\prototile'$ contain edges $\edge$ and $\edge'$ respectively such that $\edge = \tau(\edge')+\vectorv$ and $\eo_{\prototile}(\edge) = \tau(\eo_{\prototile'}(\edge'))+\vectorv$, then $\edge$ and $\edge'$ must have the same \defemph{edge breakdowns}, which means roughly that their inflated images must contain the same edges in the same order with the same orientations.  
More specifically, let $\edge_1\subset\tile_1,\ldots,\edge_k\subset\tile_k$ denote the edges in $\sub(\prototile)$ that are contained in $\infl\edge$, and let $\edge_1'\subset\tile_1',\ldots,\edge_m'\subset\tile_m'$ denote the edges in $\sub(\prototile')$ that are contained in $\infl\edge'$.  
Then $k = m$, $\edge_i = \tau(\edge_i')+\vectorv$ for all $1\leq i\leq k$, and $\eo_{\tile_i}(\edge_i) = \tau(\eo_{\tile_i'}(\edge_i))+\vectorv$ for all $1\leq i\leq k$.  
\end{enumerate}


Now let us describe the special prototile sets that we will use.  
For a natural number $n \geq 3$, let $\Triangles{n}$ denote the set of isometry classes of triangles, the angles of which are integer multiples of $\pi / n$, normalized so that they can all be inscribed in circles of the same size.  
$\Triangles{7}$ is depicted in Figure \ref{FIG:n-gon-inscribe}.  


The goal in this work is to find substitutions that use proper subsets of $\Triangles{n}$ as prototiles.  
The substitution depicted in Figure \ref{FIG:sample-substitution} uses three of the four prototiles from $\Triangles{7}$.  

Note that the general method for constructing substitutions that is described in \cite{DN:n-fold} uses a bigger set of prototiles.  
In particular, that prototile set contains two copies of each scalene triangle in $\Triangles{n}$.  
These triangles have opposite edge orientations and therefore should be considered as different from one another.  
Let us not follow this convention here; for us, any two isometric tiles will always have the same edge orientations, and hence will be copies of the same prototile.  

%% file: lengths-areas.tex
\section{Lengths and Areas}\label{SEC:lengths-and-areas}

We will need to know the areas of the triangles in $\Triangles{n}$, along with their various edge lengths.  
Denote by $\tr{k_1}{k_2}{k_3}$ the triangle with angles $k_1\pi/n,k_2\pi/n$, and $k_3\pi/n$.  

\subsection{Lengths}\label{SUBSEC:lengths}

\begin{figure}
\includegraphics[width=\textwidth]{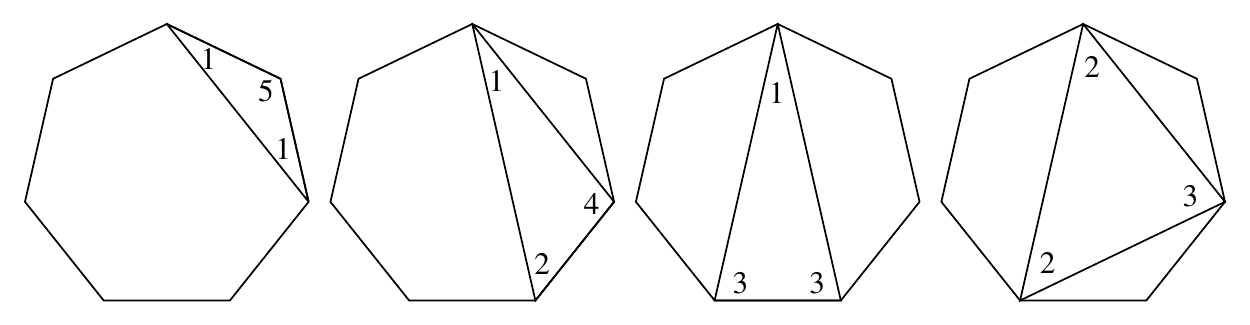}
\caption{The elements of $\Triangles{7}$ inscribed in regular heptagons.\label{FIG:n-gon-inscribe}}
\end{figure}

Since we have normalized the triangles in $\Triangles{n}$ so that they can all be inscribed in the same circle, in particular we know that their side lengths coincide with the lengths of the diagonals of a regular $n$-gon, as we have done for $n = 7$ in Figure \ref{FIG:n-gon-inscribe}.  

Let the edge lengths of this regular $n$-gon be $a_1 := 1$.  
Then the first diagonal has length $a_2 := \sin ((n-2)\pi /n)/\sin (\pi/n) = 2\cos (\pi/n)$.  
This is the length of the long edge of the triangle $\tr{1}{1}{n-2}$.  

\begin{figure}
\begin{subfigure}{\textwidth}
\centering
\includegraphics[width=1.0\textwidth,page=1]{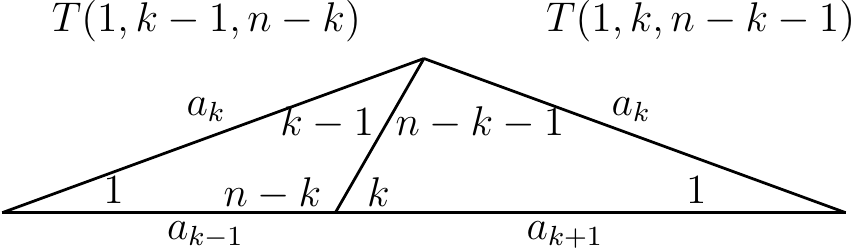}

\caption{Two narrow triangles.}
\label{FIG:similar-narrow}
\end{subfigure}
\begin{subfigure}{\textwidth}
\centering
\includegraphics[width=1.0\textwidth,page=2]{substitution-search-paper-pics.pdf}

\caption{A non-narrow triangle and an inflated narrow triangle.}
\label{FIG:similar-non-narrow}
\end{subfigure}
\caption{Similar triangles for recursion formulas.}
\label{FIG:similar}
\end{figure}

We can use similar triangles to produce a recursion relation for the length $a_k$ of the $k$th diagonal of the $n$-gon.  
In particular, the triangles $\tr{1}{k}{n-k-1}$ and $\tr{1}{k-1}{n-k}$ can be placed against one another along their short edges to produce a triangle similar to $\tr{1}{1}{n-2}$, as has been done in Figure \ref{FIG:similar-narrow}.  
The edges of this triangle are obtained from the edges of $\tr{1}{1}{n-2}$ by scaling by a factor of $a_{k}$.  
In particular, using the fact that the long edge of this new triangle has length $a_{k+1} + a_{k-1}$, we see that
\begin{align}
a_{k-1} + a_{k+1} & = a_2a_k\\
a_{k+1} & = a_2a_k - a_{k-1}.\label{EQ:recursion}
\end{align}

Taking $a_0=0$, this recursion holds for all $k\geq 1$.  
If we view $a_2$ as a variable and $a_k$ as a polynomial in $a_2$, then $a_k$ is the $k$th Chebyshev polynomial of type 2, subject to the reparametrization $a_2 = 2x$ \cite{R:chebyshev}.  
If $n$ is odd, then diagonal number $(n-1)/2$ has the same length as diagonal number $(n-1)/2 - 1$, so setting 
\begin{align}
a_{(n-1)/2} & = a_{(n-1)/2-1}\label{EQ:relation}
\end{align}
yields a polynomial equation of degree $(n-1)/2$ that $a_2$ satisfies.  

If $n$ is even, then the equation is
\begin{align*}
a_{n/2} & = a_{n/2-1}.
\end{align*}

Note that, by the symmetry of the regular $n$-gon, $a_k = a_{n-k}$.  


Now let $A$ denote the companion matrix of the minimal polynomial $q_n$ of $a_2$.  
If $b\in\Q(a_2)$, let $p_b\in\Q[x]$ denote the monic polynomial for which $b=p_b(a_2)$, and let $v(b) = (v_0(b),\ldots,v_{\Phi(n)/2}(b))$ denote the vector of coefficients of this polynomial; i.e., the vector representation of $b$ with respect to the basis $1,a_2,\ldots,a_2^{\Phi(n)/2}$ of $\Q(a_2)$.  
Note that by Equation \ref{EQ:recursion} $p_{a_i}$ has integer coefficients.  
Then $v(a_2b)^t = Av(b)^t$ and, given an inflation factor $\infl\in\Q(a_2)$, $v(\infl b)^t = p_\infl(A)v(b)^t$.  
This provides a means of representing the lengths $\lambda a_i$ as combinations of $a_0,\ldots,a_{(n-1)/2}$.  
In particular, let $L_n := [v(a_0)|v(a_2)|\cdots |v(a_{(n-1)/2}]$ denote the $\Phi(n)/2 \times (n-1)/2$ matrix, the columns of which are $v(a_i)$; then the $i$th column of the solution $X$ of the matrix equation
\begin{align}\label{EQ:lengths-composite}
L_n X & = p_\infl(A)L_n
\end{align}
expresses $\infl a_i$ as a combination of $a_0,\ldots,a_{(n-1)/2}$.  
If $n$ is prime, then $L_n$ is square---the number $(n-1)/2$ of edge lengths agrees with the degree of $q_n$---and upper triangular---the formula \ref{EQ:recursion} expresses $a_i$ as a polynomial of degree $i$ in $a_2$---and hence invertible.  
Then we can write  
\begin{align}\label{EQ:lengths-prime}
X & = L_n^{-1}p_\infl (A)L_n.
\end{align}

Let us use only inflation factors $\infl$ that are positive integer combinations of $a_0,\ldots,a_{(n-1)/2}$.  
The reason for doing so is that, if $\tile\in\Triangles{n}$ contains an angle of measure $\pi/n$, then its shortest edge has length $a_0=1$, so the shortest edge of $\infl \tile$ will have length $\infl$, which must therefore be a sum of prototile edge lengths.  
This is not strictly necessary if we choose for our set of prototiles a subset of $\Triangles{n}$ that contains no triangle with an edge of length $1$, but empirically such prototile sets do not work well, and we will soon focus our attention on a special set of prototiles that contains two triangles with minimal edge lengths---see Section \ref{SUBSEC:special}.  

\subsection{Areas}\label{SUBSEC:areas}

Now we can calculate the areas of the triangles in $\Triangles{n}$.  
Given a triangle $\tile$, let $\area{\tile}$ denote its area.  
Then we can express the areas of all the triangles as elements of $\Q(a_2)\cdot\area{\tr{1}{1}{n-2}}$, where $a_2=2\cos(\pi/n)$ as described in Section \ref{SUBSEC:lengths}.  
Let us call a triangle in $\Triangles{n}$ a \defemph{narrow triangle} if it has an angle of $\pi/n$.  
Then we first calculate the areas of the narrow triangles recursively using the same similar triangles that we used to express $a_k$ in terms of $a_2$.  

In particular, $\tr{1}{k}{n-k-1}$ and $\tr{1}{k-1}{n-k}$ fit together to form a triangle similar to $\tr{1}{1}{n-2}$, but inflated by a factor of $a_{k}$ (see Figure \ref{FIG:similar-narrow}).  
Therefore 
\begin{align}
\area{\tr{1}{k}{n-k-1}} & = a_{k}^2\cdot\area{\tr{1}{1}{n-2}} - \area{\tr{1}{k-1}{n-k}},
\end{align}
so, recursively, we obtain formulas expressing the areas of the narrow triangles in terms of $\area{\tr{1}{1}{n-2}}$.  
Note that each $a_k$ can be written as an integer polynomial in $a_2$, so this recursion formula expresses $\area{\tr{1}{k}{n-k-1}}$ as an integer polynomial in $a_2$.  

To calculate the areas of the non-narrow triangles, note that every non-narrow triangle fits together with an inflated narrow triangle to form another inflated narrow triangle, as in Figure \ref{FIG:similar-non-narrow}.  
In particular, if $k<l<n-k-l$, then 
\begin{align}
\area{\tr{k}{l}{n-k-l}} & = a_k^2\cdot\area{\tr{1}{l}{n-l-1}} - a_l^2\cdot\area{\tr{1}{k-1}{n-k}},
\end{align}
and, since $\area{\tr{1}{l}{n-1-l}}$ and $\area{\tr{1}{k-1}{n-k}}$ are multiples of $\area{\tr{1}{1}{n-2}}$ by integer polynomials in $a_2$, so is $\area{\tr{k}{l}{n-k-l}}$.  
Note that permuting the order of the angles $k$, $l$, and $n-k-l$ gives us up to three different formulas for this area, any one of which will work.  

The main reason for computing triangle areas is to determine the substitution matrix of the substitution $\sub$ that we are trying to find.  
Choose a set of prototiles $\PT\subset \Triangles{n}$, the areas of which form a $\Q$-basis for $\Q(a_2)\cdot \area{\tr{1}{1}{n-2}}$.  
Enumerate the elements of this set $\tile_1,\ldots,\tile_m$, and let $A_k := \area{\tile_k}/\area{\tr{1}{1}{n-2}}$.  
Let $B_\PT$ denote the square matrix, the columns of which are the vectors $v(A_k)$ (defined in Section \ref{SUBSEC:lengths}); that is, $B_\PT = [v(A_1)|\cdots |v(A_m)]$.  
Our choice of $\PT$ having areas that are a $\Q$-basis for $\Q(a_2)\cdot\area{\tr{1}{1}{n-2}}$ means that $B_\PT$ is invertible.  

Then, if there exists a substitution $\sub$ on $\PT$ with inflation factor $\infl$, it will have substitution matrix
\begin{align}\label{EQ:substitution-matrix}
M & := B_\PT^{-1}p_\infl(A)^2B_\PT.
\end{align}

Sometimes an inspection of the matrix $M$ is enough to prove the non-existence of a substitution on $\PT$ with factor $\infl$.  
For instance, if any of the entries of $M$ are not integers or are negative, then no such rule can exist---although \cite{F:thesis} describes a way of modifying $\PT$ to address the problem of a negative entry.  

If the entries of $M$ are all positive integers, then we search the combination of $\PT$ and $\infl$ for substitution rules.  

%% file: algorithm.tex
\section{The Program}\label{SEC:algorithm}

Let us now describe the program used to search for substitution rules.  
The program is written in java, although the discussion presented here contains no java-specific details.  

\subsection{Overview}\label{SUBSEC:overview}

Before describing the search algorithm in any detail, let us give an overview of how the program works.  

The order of symmetry $n$ is assumed to be fixed from the start.  
The program takes three ingredients as input: a set $\PT\subset\Triangles{n}$ of prototiles, an inflation factor $\infl$ that is a non-negative integer combination of $a_1,a_2,\ldots,$ $a_{(n-1)/2}$, and an inflated prototile $\infl\tile_0$, where $\tile_0\in\PT$.  
Then it tries to fill $\infl\tile_0$ with tiles congruent to prototiles from $\PT$ in such a way that the tiles overlap at most in their boundaries, they meet edge-to-edge if at all, and their union is all of $\infl\tile_0$.  
Let us refer to each such application of the program as a \defemph{search} of the triple $(\PT,\infl,\tile_0)$.  
The output of a search is a (possibly empty) set of patches from $\PatchSet(\Triangles{n})$, each of which has support $\infl\tile_0$.  
Let us call such a patch a \defemph{result}.  
A search will be considered successful if it returns at least one result.  

In order to obtain a substitution rule on $\PT$ with factor $\infl$, we must run a search on $(\PT, \infl, \tile)$ for each prototile $\tile\in\PT$, and the search must be successful for each one.  
Then we can define a substitution rule $\sub$ as follows: for each $\tile\in\PT$, let $\sub(\tile)$ be any one of the patches found in the search on $(\PT,\infl,\tile)$.  
In order for the substitution rule so obtained to be edge-to-edge, it must satisfy some additional conditions on its edge orientations and edge breakdowns, as described in Section \ref{SUBSEC:objective}.  
Isometry equivariance (condition (1)) and matching (condition (2)) are built into the search---so, in fact, the program will not find any patches that do not satisfy these conditions, although that restriction can be turned off.  

Preservation under $\sub$ (condition (3)) is checked after the completion of all the searches for a given $\PT$ and $\infl$.  
More specifically, we assemble all of the results together and look for a tuple of results, one for each tile, that can be used to define a substitution that satisfies condition (3).  
We also check now that the tuple of results satisfies condition (1) as a whole, because the program only checks the patches individually during the search.  

If such a tuple of results exists, then the substitution that it defines is edge-to-edge, and we have found what we set out to find for $\PT$ and $\infl$.  

\subsection{A Special Prototile Set}\label{SUBSEC:special}

Of all subsets $\PT\subset \Triangles{n}$ for odd $n$, one works better than the others.  
This is the set consisting of all triangles that contain at least one edge of maximum length, that is, all triangles with at least one angle of measure $n/2 \pm 1/2$.  
There are $(n-1)/2$ such triangles, and they can be combined along their short sides to produce triangles similar to the $(n-1)/2$ isosceles triangles in $\Triangles{n}$.

\subsection{Representation of Points}\label{SUBSEC:points}

If a patch in $\PatchSet(\Triangles{n})$ contains a tile with a vertex at the origin, then all tile vertices in that patch must lie in the $\Z$-module $\Vertices$ generated by all of the vectors
\begin{align}\label{EQ:generators}
\{ a_k (\cos (i\pi/n),\sin(i\pi/n)) \ | \ 1\leq k \leq (n-1)/2, \ 0\leq i < 2n\}.
\end{align}
If each $\tile\in\PT$ has a vertex at the origin, then so does each $\infl\tile$, and so we can restrict our attention to such patches.  
Therefore the program need not deal with arbitrary points in $\R^2$, but only points in $\Vertices$.  

Internally, it represents these points as elements of $\Z^{n-1}$.  

\begin{figure}
\includegraphics[width=1.0\textwidth,page=3]{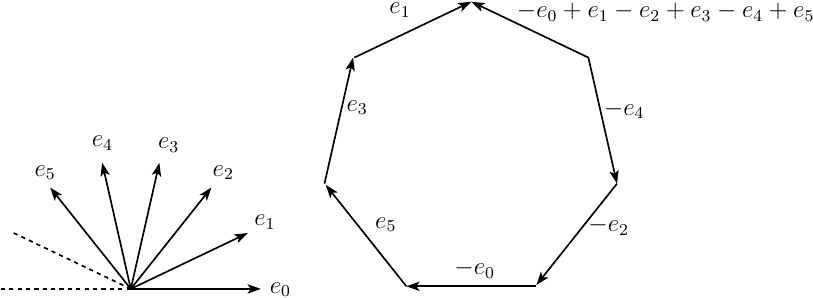}

\caption{The images of the standard basis vectors of $\Z^6$, arranged to form the edges of a regular heptagon.}
\label{FIG:vectors}
\end{figure}

If $n$ is prime, then $\Z^{n-1}$ is isomorphic to $\Vertices$.  
Specifically, multiplying on the left by the matrix
\begin{align}\label{EQ:matrix}
\left[ \begin{array}{cccc}
\cos \frac{0\pi}{n} & \cos \frac{1\pi}{n} & \cdots & \cos \frac{(n-2)\pi}{n} \\
\sin \frac{0\pi}{n} & \sin \frac{1\pi}{n} & \cdots & \sin \frac{(n-2)\pi}{n} 
\end{array}\right]
\end{align}
takes the $i$th standard basis vector of $\Z^{n-1}$ to the vector $(\cos (i\pi/n),\sin (i\pi/$ $n))$ $\in\Vertices$, which is an element of the set \ref{EQ:generators} as $a_1=1$.  
To see that this map is onto, note that the images of the standard basis vectors are the direction vectors of $n-1$ of the sides of a regular $n$-gon.  
This is depicted for the case $n=7$ in Figure \ref{FIG:vectors}.  
Then it is not hard to see that vectors with all different lengths $a_1,\ldots,a_{(n-1)/2}$ can be obtained as sums of the images of the standard basis vectors---for instance, 
\begin{align*}
a_2(\cos(1\pi/n),\sin(1\pi/n)) & = (\cos(0\pi/n),\sin(0\pi/n)) \\ & \quad + (\cos(2\pi/n),\sin(2\pi/n)).
\end{align*}
If $n$ is not prime then this map $\Z^{n-1}\to \Vertices$ is onto but not one-to-one, as there are multiple non-trivial relations between the direction vectors of the edges of a regular $n$-gon.  

\subsection{The Search}\label{SUBSEC:search}


The search works by recursion.  
It has two main objects that it updates with each recursion step, and it uses these objects to decide whether or not to terminate the recursion.  
The first object is a patch $\Patch\in\PatchSet(\PT)$, the support of which is contained in $\infl\tile_0$.  
The second is a list of non-negative integer multiplicities of prototiles.  
These multiplicities indicate how many prototiles of each type remain to be placed in order to obtain a patch with support equal to $\infl\tile_0$; it is initialised using the column of the matrix \ref{EQ:substitution-matrix} that corresponds to $\infl\tile_0$.  
We say that this list contains a tile $\tile$ if the multiplicity of that tile is greater than $0$.  
At each step in the recursion, the multiplicity of one prototile is decremented in this list and the corresponding tile is pushed on $\Patch$, which is in fact a stack.  
Then the recursion is finished when all of the multiplicities are $0$.  


The main recursive procedure is called solve.  
Simplified pseudocode for this procedure appears in Algorithm \ref{ALG:pseudocode}.  

\begin{algorithm}
\caption{The solve procedure}\label{ALG:pseudocode}
\begin{algorithmic}[1]
\State{$\Patch$: a patch with support contained in $\infl\tile_0$}
\State{$\tile$: the prototile we are currently trying to place in $\Patch$}
\State{$l$: a list of multiplicities telling us how many prototiles of each type remain to be placed}
\Procedure{solve}{}
\Repeat
\If{$l$ is empty} 
\State{store a deep copy of $\Patch$ somewhere}
\Comment{$\Patch$ is a result}
\State{\textbf{break}}
\EndIf
\If{$l$ contains $\tile$}
\State{$\tile'$ $\gets$ place $\tile$}\label{line:place}
\If{$\tile'$ is compatible with $\Patch$}\label{line:compatible}
\State{push $\tile'$ on $\Patch$}\label{line:forward}
\State{$\tile$ $\gets$ first prototile}\label{line:reset}
\State{solve}\label{line:recursion}
\Comment{recursion}
\State{$\tile$ $\gets$ pop last tile from $\Patch$, get associated prototile}\label{line:back}
\EndIf
\EndIf
\State{$\tile$ $\gets$ prototile after $\tile$}\label{line:next-prototile}
\Until{back to first prototile}\label{line:first-prototile}
\EndProcedure
\end{algorithmic}
\end{algorithm}

Several of the statements in Algorithm \ref{ALG:pseudocode} require further explanation.  
The place procedure on line \ref{line:place} takes as input the prototile $\tile$ and returns a tile $\tile'$ that is congruent to $\tile$.  
Then the compatible procedure on line \ref{line:compatible} returns true if $\Patch \cup \{\tile'\}$ is a patch with support contained in $\infl\tile_0$.  
(In fact, compatible also requires $\Patch\cup\{\tile'\}$ to satisfy edge orientation conditions (1) and (2)---see Sections \ref{SUBSEC:objective} and \ref{SUBSEC:overview}.)  


There are many conceivable ways of placing a copy of the prototile $\tile$ in $\Patch$, but we only consider a restricted range of possibilities.  
The edges of all of the tiles in $\Patch$ are divided into two lists, called open edges and closed edges.  
Closed edges are edges that are contained in two different tiles, or in one tile and in the boundary of $\infl\tile_0$.  
Open edges are edges that are contained in only one tile, and not in the boundary of $\infl\tile_0$.  
In other words, open edges are edges against which another tile must be placed in order to obtain a completed patch.  

The open edges and closed edges are stored in stacks, so in particular they have an order, with the open edge(s) coming from the most recently-placed tile being listed last.  
Then the rule is that we place the prototile $\tile$ against the last open edge, if this is possible ($\tile$ might not have an edge of the right length).  
This rule has two main implications.  

The first implication is that any newly-placed tile $\tile'$ has at least one edge in common with a tile that has already been placed---this edge becomes a closed edge.  
Therefore a newly-placed tile can contribute at most two new open edges.  
If it contributes one new open edge, or no new open edges, then there is no ambiguity as to which new open edge is listed last.  
But if it contributes two new open edges, then we have a systematic way of determining the order in which they are added to the stack; namely, they are pushed in the order in which they appear in a counterclockwise traversal of the edges of $\tile'$, starting with the closed edge.  

\begin{figure}
\includegraphics[width=1.0\textwidth,page=4]{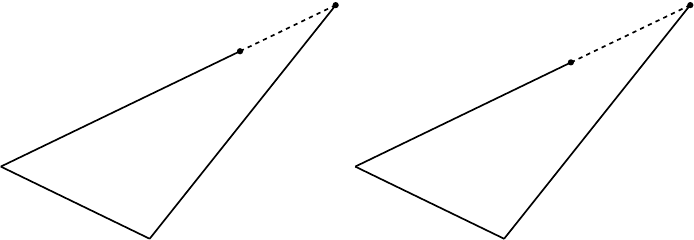}

\caption{Two copies of an inflated prototile with starter edges of different lengths shown as dashed lines.}
\label{FIG:starters}
\end{figure}

The second implication is that, in order to call the place procedure, we need to have at least one open edge, even at the beginning when no tiles have been placed yet.  
This is a problem that we fix by placing an open edge, called a starter, on part of one of the edges of $\infl\tile_0$ (see Figure \ref{FIG:starters}).  
Then the first tile must be placed against the starter.  

We need to choose a length for the starter edge, and some choices of edge length might preclude certain results.  
We deal with this by running several searches in parallel---one for each choice of starter edge length.  
This turns out to be a convenient method to divide up the work for multithreading---see Section \ref{SUBSEC:multithreading}.  
To produce a list of possible starter edge lengths, one option is to use all edge lengths $a_1,\ldots,a_{(n-1)/2}$.  
An even better option is to select only those lengths for which the corresponding entry in the length substitution matrix $X$ from Equation \ref{EQ:lengths-prime} (Equation \ref{EQ:lengths-composite} if $n$ is not prime) is positive.  

The pseudocode in Algorithm \ref{ALG:pseudocode} has been simplified a bit for clarity, but we should note now that there might be several ways to place the prototile $\tile$ against the last open edge.  
If $\tile$ has an edge of the same length as the last open edge, then there are two possibilities: we can place $\tile$ against the last open edge using a non-orientation preserving isometry---i.e., a reflection---or using an orientation preserving isometry.  
If $\tile$ has two edges of that length, then there are four ways of placing it against the last open edge.  
(The case in which $\tile$ is equilateral has not come up yet, because the program only works on prime $n$ now, so in particular not on any $n$ divisible by $3$.)  

To keep track of these possible placements, we store not only the current prototile $\tile$, but also two booleans $\flip$ and $\second$ that tell us respectively if we are trying to place a reflected copy of $\tile$ and if we are trying to place $\tile$ using the second of two edges of the same length.  
Therefore on line \ref{line:place}, we do not place $\tile$, but rather we place the triple $(\tile,\flip,\second)$.  
On line \ref{line:next-prototile}, we get the next triple $(\tile,\flip,\second)$, not just the next prototile $\tile$, and on lines \ref{line:reset} and \ref{line:first-prototile} we refer to the first triple $(\tile_1,false,false)$, not just the first prototile $\tile_1$.  
On line \ref{line:back}, we pop the last tile $\tile'$ from $\Patch$ and, using information from $\tile'$ and $\Patch$ we infer which triple $(\tile,\flip,\second)$ was input to the place procedure to produce $\tile'$.  


This last point merits further discussion.  
The purpose of line \ref{line:back} is to restore the states of $(\tile,\flip,\second)$, and $\Patch$ to what they were on line \ref{line:forward}.  
In principle we could store these states in memory, but note that we call the solve procedure recursively on line \ref{line:recursion}; this means that these states could change many times between lines \ref{line:forward} and \ref{line:back}.  
In practice it becomes unwieldy to store all of these changes---sometimes even crashing the program ---so we store none of them.  
Instead, when we remove $\tile'$ from $\Patch$ on line \ref{line:back}, we use information about $\tile'$ to restore everything to its state from line \ref{line:forward}.  

In particular, we remove any open edges coming from edges of $\tile'$, and if any closed edges come from edges of $\tile'$, then they are made open again.  
We also determine the prototile $\tile$ to which $\tile'$ is congruent.  
If $\tile'$ is a reflected copy of $\tile$, then we set $\flip$ to $true$, and if $\tile'$ has its second of two equal edges against the last open edge, then we set $\second$ to $true$.  
These calculations allow us to ascend and descend the recursion tree without having to use up the system memory storing states from previous levels of recursion.  

\subsection{Multithreading}\label{SUBSEC:multithreading}

The amount of time required to search an inflated prototile $\infl\tile_0$ scales up rapidly with its area, and it is worthwhile to take advantage of the multithreading capabilites of the java language to reduce this time.  
The particular algorithm design that makes it possible to do this might be of independent interest, so let us describe it here.  

The program holds a list of searches at varying stages of completion; at the same time, there are many threads, each working on a different search.  
When a thread finishes a search, it draws a new one from the list.  

These searches are independent from one another, in the sense that no two can ever return the same result.  
This is because they all contain different configurations of open edges and tiles that have already been placed in their patches.  
So the initial searches are not actually empty; rather, they contain starter edges (see Section \ref{SUBSEC:search}), which differ from one another in their lengths.  
Thus we can have up to $(n-1)/2$ initial searches.  

We may have many more threads than initial searches, and, empirically, most of the initial searches seem to finish quickly with no results.  
Therefore we need a way to split up the few searches that remain in order to take advantage of the extra threads available.  
We do this by giving each search a kill switch that we trigger when the number of searches remaining in the list drops too low.  

When the kill switch is triggered, the search algorithm changes: instead of making another solve call at line \ref{line:recursion} of Algorithm \ref{ALG:pseudocode}, the search makes a deep copy of itself in its current state, and places this copy in the list of searches.  
A given search may have called solve recursively many times already when the kill switch is triggered, so it will pass line \ref{line:recursion} many times before finishing, adding many searches to the list.  
These new searches represent work that the original search would have done had the kill switch not been triggered.  

Figuratively, the behavior of the searches can be described as follows.  
Under ordinary conditions, a search climbs the recursion tree and searches all of the branches to their ends as they are encountered.  
After the kill switch has been triggered, the search no longer climbs up new branches that it encounters; instead, whenever it encounters a new branch on the way down the tree, it cuts it off at the base and throws it on a pile of branches to be searched.  
Then it reaches the bottom of the recursion tree and finishes.  

%% file: results.tex
\section{Results}\label{SEC:results}

\subsection{Fivefold Symmetry}\label{SUBSEC:fivefold}

To find substitutions using $\Triangles{5}$ it is not necessary to use a computer, but it certainly helps.  
For these prototiles, the smallest inflation factor is $a_2$, the golden ratio.  
Substitutions with this inflation factor have been described elsewhere (see \cite{DN:n-fold} and \cite{G:fivefold}); these can be combined as multi-substitutions, but not as random substitutions---although see \cite{G:fivefold} for a description of a slightly different randomization strategy, called tile rearrangement.  

\begin{figure}
\includegraphics[width=1.0\textwidth]{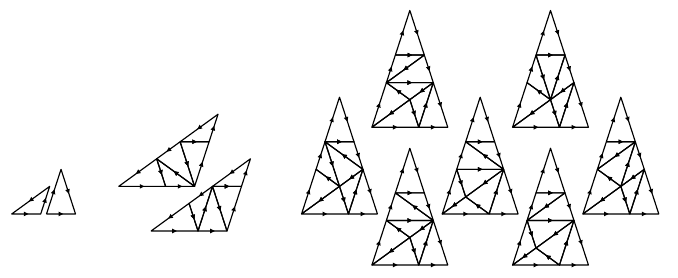}
\caption{Results for $\Triangles{5}$, $\infl = 1+a_2$.}
\label{FIG:fivefold}
\end{figure}

Already at the next inflation factor, $a_2^2 = 1+a_2$, there are many substitutions to be found.  
Figure \ref{FIG:fivefold} depicts a family of these---two for $\tr{1}{1}{3}$ and seven for $\tr{1}{2}{2}$---that can be used to produce random substitutions.  
The program also finds substitutions with different edge breakdowns that can be combined with these ones as multi-substitutions.

\subsection{Sevenfold Symmetry}\label{SUBSEC:sevenfold}

We have found many new substitution rules with sevenfold rotational symmetry using three different inflation factors: $\infl_1=1+a_2$, $\infl_2=a_2+a_3$, and $\infl_3=1+a_2+a_3$.  
The second and third of these are both PV numbers, while the first is not.  
A PV inflation factor is a necessary condition for the existence of any non-trivial eigenvalue of the dynamical system, that is, for any non-trivial discrete part in the spectrum \cite{S:pisot}.
All three inflation factors are units in the ring of integers of the number field generated by $a_2$, which is necessary for the module generated by the eigenvalues to be finitely generated.  

Of course, it is possible to search even larger inflation factors, and presumably to find more substitution rules, at the cost of more computing time and much greater memory usage.  

The $\tr{1}{3}{3}$ and $\tr{2}{2}{3}$ triangles are both isosceles, so each one has two ways of placing it in any given position: an orientation preserving placement and an orientation reversing placement.  
If we choose one of these isosceles triangles and reflect each of its instances in a result $\Patch$ across its axis of symmetry, then this modified patch is still a result.  
We have processed the results obtained from the program to remove this redundancy and to select only those groups of results that have compatible edge breakdowns in the sense of Section \ref{SUBSEC:objective}, and that therefore can be used to produce tilings.  

\begin{figure}
\begin{subfigure}{\textwidth}
\includegraphics[width=1.0\textwidth]{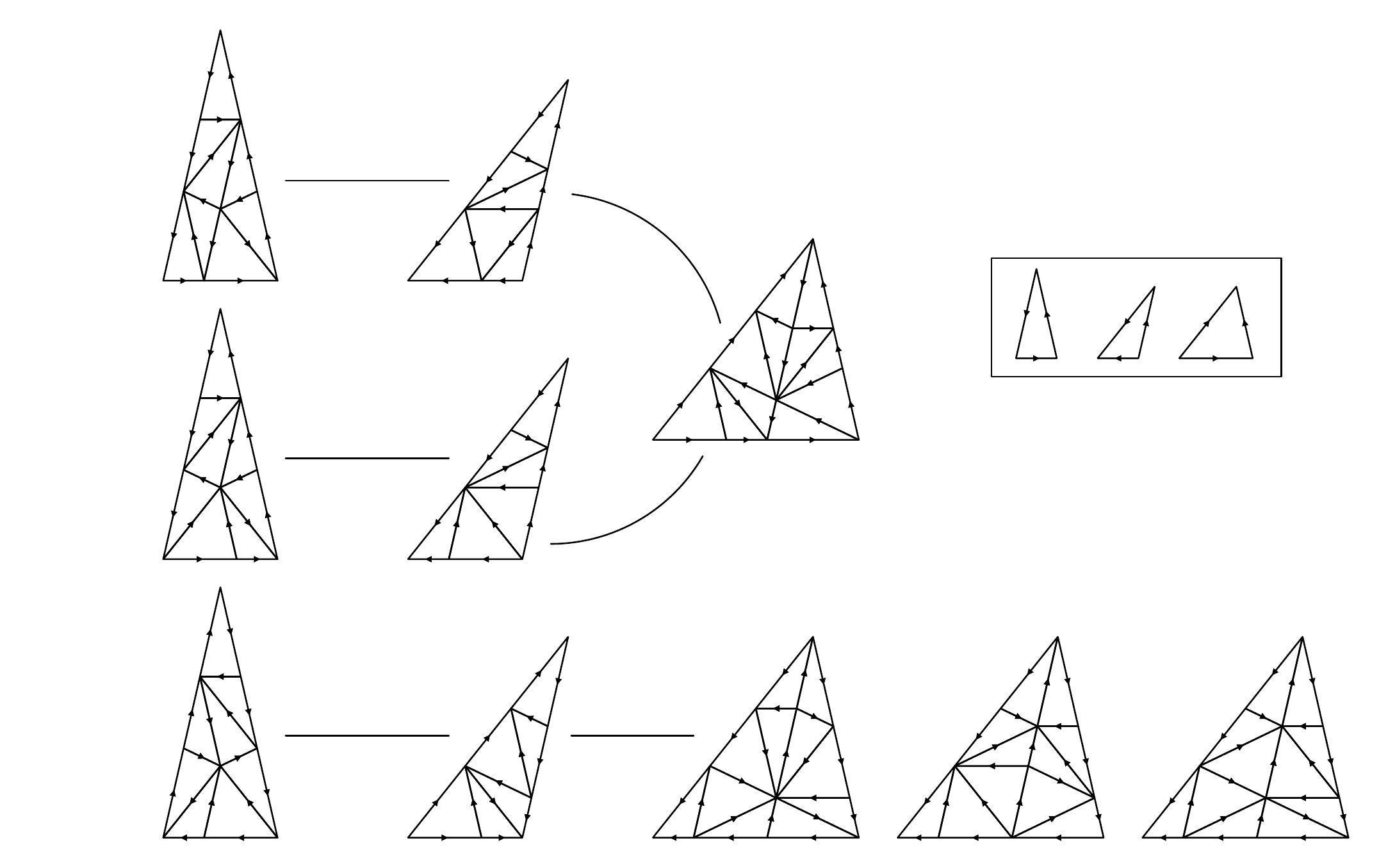}
\caption{The standard prototiles.  These are compatible with results for $\infl_2$ and $\infl_3$.}
\label{FIG:small-sevens-standard}
\end{subfigure}

\begin{subfigure}{\textwidth}
\includegraphics[width=1.0\textwidth]{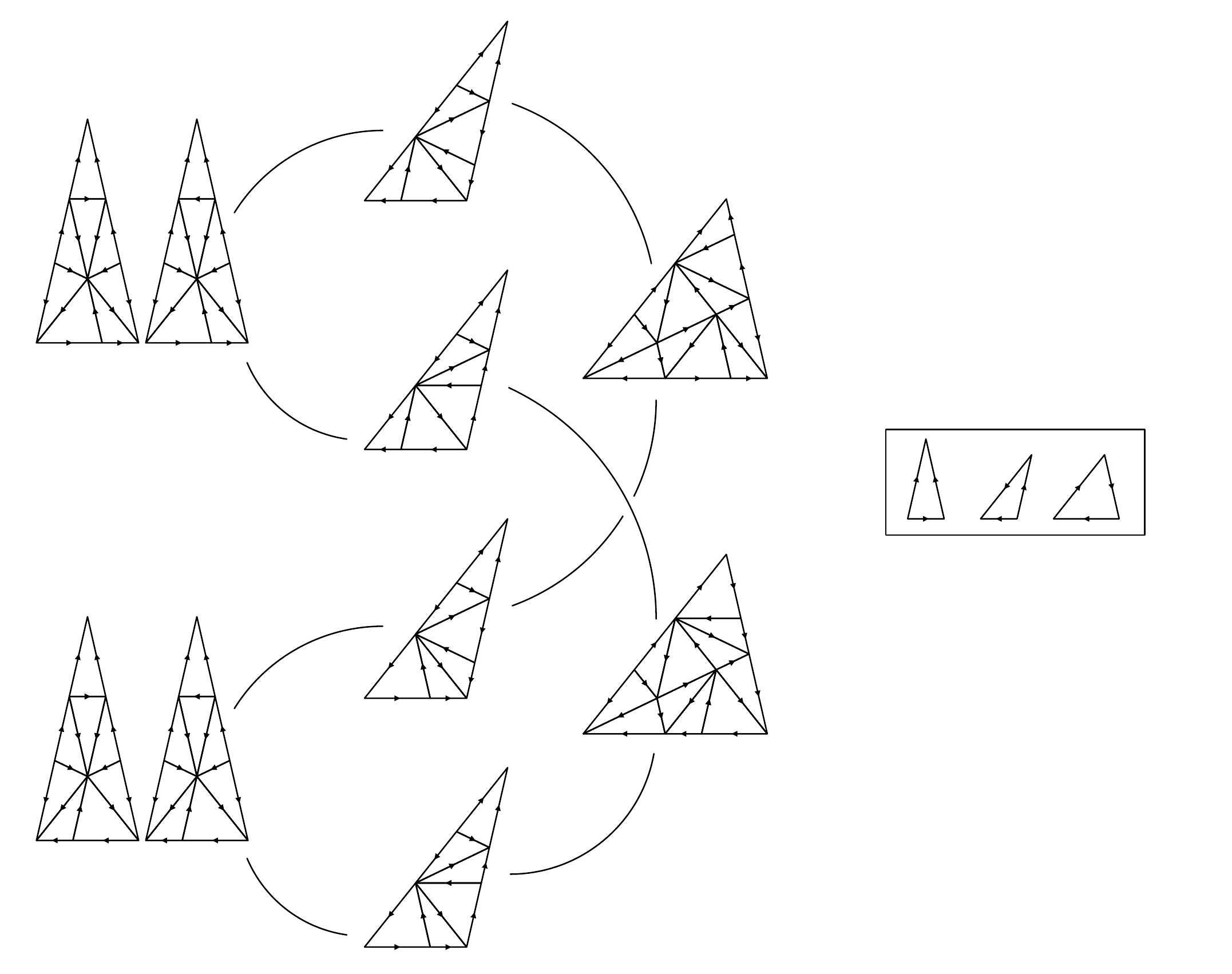}
\caption{The non-standard prototiles.  The Danzer-Nischke example \cite{DN:n-fold} uses the result for $\tr{1}{2}{4}$ depicted at the top.}
\label{FIG:small-sevens-different}
\end{subfigure}
\caption{All results for inflation factor $\infl_1$.  They fall into two subsets with different edge orientations on the prototiles, which appear in boxes on the right.  Arcs between patches indicate compatibility of edge breakdowns.}
\label{FIG:small-sevens}
\end{figure}

The results for $\infl_1$ appear in Figure \ref{FIG:small-sevens}.  
They fall into two classes, depicted in Figures \ref{FIG:small-sevens-standard} and \ref{FIG:small-sevens-different}, according to the edge orientations of the prototiles, which are enclosed in boxes on the right sides of the figures.  
These should really be considered as two different sets of prototiles, because substitution rules coming from the two different classes of results cannot be combined with one another even as multi-substitutions.  
On the other hand, two or more substitution rules coming from the same class of results can be combined as multi-substitutions.  

An arc between two results in Figure \ref{FIG:small-sevens} indicates that they have the same edge breakdowns.  
Two or more substitution rules coming from results with the same edge breakdowns can be combined as random substitutions.  
So, for instance, we can make three different substitution rules by choosing, from the bottom row of Figure \ref{FIG:small-sevens-standard}, any one of the three results for $\tr{2}{2}{3}$, along with the results for $\tr{1}{3}{3}$ and $\tr{1}{2}{4}$; these three substitutions can be combined as random substitutions.  
This involves creating tilings by repeatedly substituting a patch, and for each $\tr{2}{2}{3}$ tile appearing in that patch, choosing at random which of the three rules to apply to it.  

\begin{figure}
\begin{center}
\hbox{%
\hspace{-25pt}
\includegraphics[width=1.3\textwidth]{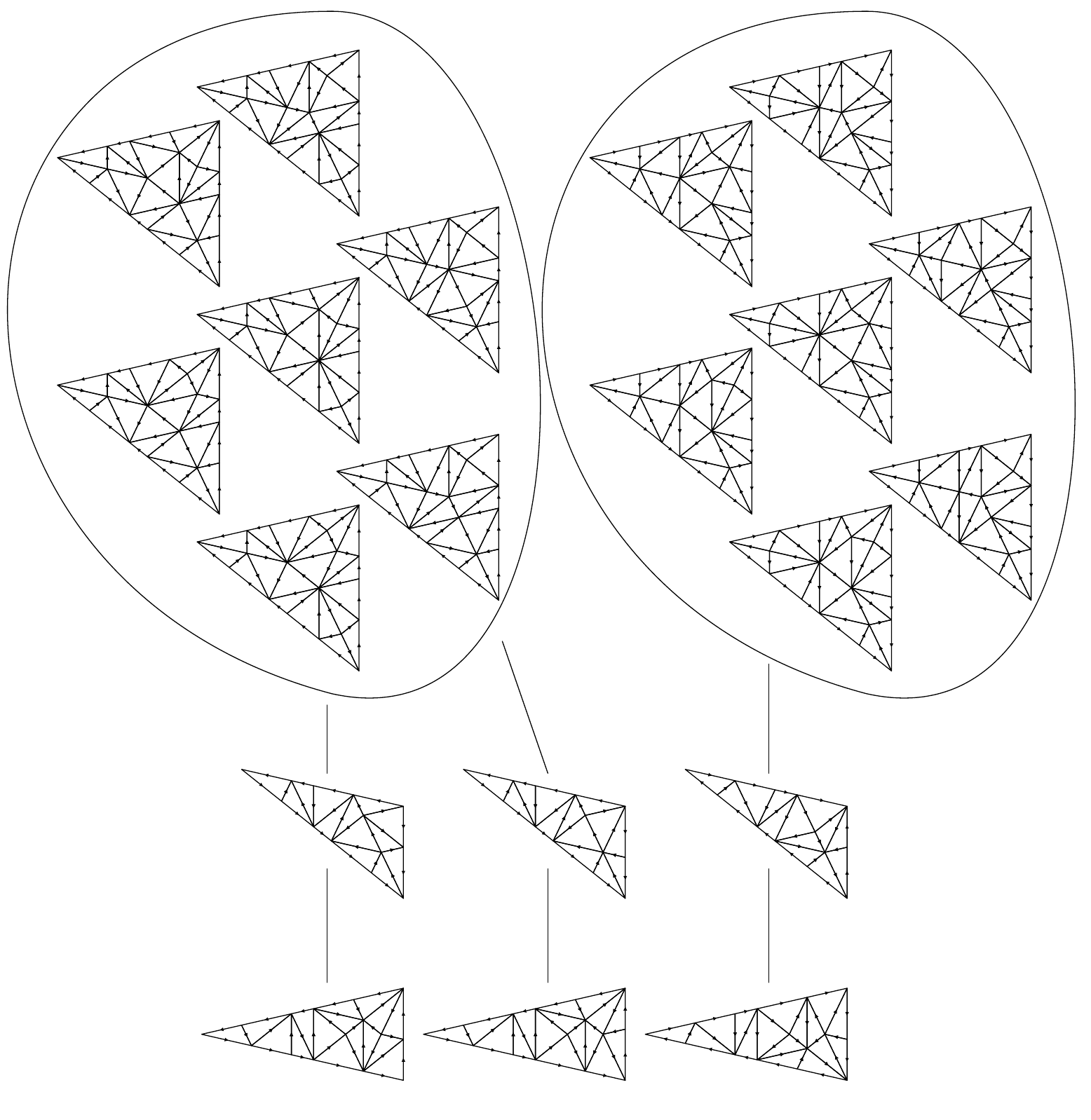}
}%
\caption{The results for $\infl_2$.  Arcs between results or groups of results indicate compatibility of edge breakdowns.  The two groups of results for $\tr{2}{2}{3}$ are reflections of one another.}
\label{FIG:medium-sevens}
\end{center}
\end{figure}

The results for $\infl_2$ appear in Figure \ref{FIG:medium-sevens}.  
The substitutions built from the results in Figures \ref{FIG:small-sevens-standard} and \ref{FIG:medium-sevens} use the same prototiles, and can be combined with one another to produce multi-substitutions.  
Any substitution built from results in Figure \ref{FIG:medium-sevens} with the same edge breakdowns and edge orientations can be combined with one another to produce random substitutions.  
Of course, a substitution built from results in Figure \ref{FIG:medium-sevens} cannot be combined with a substitution built from results in Figure \ref{FIG:small-sevens} to produce a random substitution because they have different inflation factors.

\begin{figure}
\begin{center}
\includegraphics[width=\textwidth]{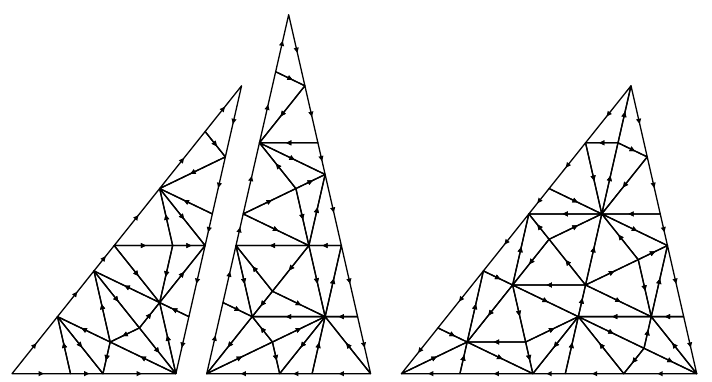}
\caption{Results for $\infl_3$.}
\label{FIG:large-sevens}
\end{center}
\end{figure}

The results for $\infl_3$ are too numerous to depict here, but a set of three of them, one for each prototile, appears in Figure \ref{FIG:large-sevens}.  
These three all have compatible edge breakdowns and edge orientations, in the sense of \ref{SUBSEC:objective}.  
For these edge breakdowns and orientations there are 3 compatible $\tr{1}{2}{4}$-results, 5 compatible $\tr{1}{3}{3}$-results, and 36 compatible $\tr{2}{2}{3}$-results, all of which can be combined with each other as random substitutions.  

Moreover, this is one out of ninety edge breakdown/edge orientation combinations; the other eighty-nine also have many results associated to them (some of which are repeated, since, for example, the $\tr{2}{2}{3}$ result appearing in Figure \ref{FIG:large-sevens} is compatible with any $\tr{1}{2}{4}$ result having the same edge breakdown on its medium and long edges, regardless of the edge breakdown of the short edge).  
A substitution built using results from any of these ninety edge breakdown/edge orientation combinations can be combined with any other such substitution, or any substitution coming from Figure \ref{FIG:small-sevens-standard} or Figure \ref{FIG:medium-sevens}.  

\subsection{Observations}\label{SUBSEC:observations}

Let us now make a few miscellaneous observations on the results presented in Sections \ref{SUBSEC:fivefold} and \ref{SUBSEC:sevenfold}.  

\begin{enumerate}
\item  The non-standard prototiles in Figure \ref{FIG:small-sevens-different} appear to be a phenomenon unique to the factor $\infl_1$, in that there are no results for the factors $\infl_2$ or $\infl_3$ with those prototiles.  
\item  The prototiles in Figures \ref{FIG:fivefold}, \ref{FIG:small-sevens-standard}, \ref{FIG:medium-sevens}, and \ref{FIG:large-sevens} have the property that, if two edges meet at a vertex at an angle that is an even integer multiple of $\pi/n$, then they are either both oriented towards that vertex or both oriented away from it; if the angle is an odd integer multiple of $\pi/n$, then their orientations are opposite.  
This has the consequence that, in any tiling made with these prototiles, all edges with the same angle will have the same orientation.  
It also means that these tilings are only $n$-fold symmetric, not $2n$-fold symmetric, as is often the case when one tries to construct something $n$-fold symmetric.  

This is also true of substitutions arising from the general method described in \cite{DN:n-fold}, although it is not true of the special substitution in that work, nor of the substitutions in Figure \ref{FIG:small-sevens-different}.  
\item  Some of the results in Figures \ref{FIG:small-sevens-standard}, \ref{FIG:medium-sevens}, and \ref{FIG:large-sevens} reverse the directions of the tile edges in the super-tile edges, while others preserve them.  
\item  The program produces many more results than are recorded here.  
Some of these results are obtained by reflecting all instances of an isosceles tile, as described in Section \ref{SUBSEC:sevenfold}, and so in that sense are redundant.  
Other results are only partial in the sense that they cannot be combined with results for the other prototiles to produce a substitution rule.  
For example, there exist results for $\tr{2}{2}{3}$ that are compatible with other results for $\tr{1}{2}{4}$, but not with any result for $\tr{1}{3}{3}$.  
In principle such results could still be used in random substitutions, provided that they were not applied to any instance of $\tr{2}{2}{3}$ that shared an edge with an instance of $\tr{1}{3}{3}$.  
\end{enumerate}

\subsection{A Surprising Topological Property}\label{SUBSEC:cohomology}

Some of these new examples produce substitution tiling spaces with a surprising topological property.  

These scaling factors are units in the relevant ring, which means that the induced substitution on the edge lengths (and the first cohomology of the Anderson-Putnam complex \cite{AP}) is unimodular, and so the first cohomology is finitely generated, and so is the module of eigenvalues in the dynamical spectrum (provided the factor is PV, otherwise there are no non-trivial eigenvalues).  
Despite all that, and unlike all examples found in the literature so far, (some of) these substitutions have a non-unimodular action on the second cohomology, which is thus not finitely generated.

All of the substitutions coming from Figure \ref{FIG:small-sevens-different} have this property, as does the substitution from the bottom row of Figure \ref{FIG:small-sevens-standard} that uses the leftmost result for $\tr{2}{2}{3}$.  
Some of the substutitions from Figure \ref{FIG:medium-sevens} have this property and others do not.  

It is also easy to find examples with this property using the fivefold prototiles $\Triangles{5}$.